\documentclass[12pt]{article}
\usepackage{amsmath,amssymb,theorem}
\usepackage{graphicx}
\usepackage{epstopdf}
\usepackage{psfrag}
\usepackage[ruled, lined, linesnumbered]{algorithm2e}
\newtheorem{thm}{Theorem}[section]

\newtheorem{lem}[thm]{Lemma}
\theorembodyfont{\rmfamily}

\begin{document}
\title{A Note on Weighted Rooted Trees }
\author{
Zi-Xia Song\thanks{Corresponding Author: Zixia.Song@ucf.edu}, Talon Ward and Alexander York\\
Department  of Mathematics\\
 University of Central Florida\\
Orlando, FL 32816, USA\\
}
\maketitle
\begin{abstract}

Let $T$ be a tree  rooted at $r$. Two vertices of $T$ are {\it related} if one is a descendant of the other; otherwise, they are {\it unrelated}. Two subsets $A$ and $B$ of $V(T)$ are {\it unrelated} if,  for any $a\in A$ and $b\in B$, $a$ and $b$  are unrelated. Let $\omega$ be a nonnegative weight function defined on  $V(T)$ with $\sum_{v\in V(T)}\omega(v)=1$. 
In this note,  we prove that either 
 there is an   $(r, u)$-path $P$ with   $\sum_{v\in V(P)}\omega(v)\ge \frac13$ for some $u\in V(T)$, or 
there exist  unrelated sets $A, B\subseteq V(T)$ such that   $\sum_{a\in A }\omega(a)\ge \frac13$ and  $\sum_{b\in B }\omega(b)\ge \frac13$. The bound   $\frac13$ is tight. This answers a question posed in a very recent paper of Bonamy, Bousquet and Thomass\'e.

\end{abstract}

\section{Introduction}




Let $T$ be a tree  rooted at $r$. Let $x\in V(T)$.  A {\it descendant} of   $x$ is any vertex $y$ such that $x\in V(P)$, where  $P$ is  the unique $(r,y)$-path in $T$.    The {\it parent} of   $x$ is the vertex $y$ such that  $y$ immediately precedes $x$ on  the unique $(r,x)$-path in $T$.  Two vertices of $T$ are {\it related} if one is a descendant of the other; otherwise,  they are {\it unrelated}. Two subsets $A$ and $B$ of $V(T)$ are {\it unrelated} if,  for any  $a\in A$ and $b\in B$, $a$ and $b$  are unrelated. Note that if $A$ and $B$ are unrelated, then $A\cap B=\emptyset$.  Let $G$ be a graph and let $\omega$ be a nonnegative weight function defined on  $V(G)$. For any $A\subseteq V(G)$ and  any  subgraph $H$ of $G$,  define  $\omega(A):=\sum_{a\in A}\omega(a)$ and $\omega(H)=\omega(V(H))$.  In their proof of the main result in~\cite{BBT}, Bonamy, Bousquet and Thomass\'e made use of the following  lemma. 

\begin{lem}\label{L}

Let $T$ be a tree  rooted at $r$ and let $\omega$ be a nonnegative weight function defined on  $V(T)$ with $\omega(T)=1$.  Then there is an $(r, u)$-path $P$ with $\omega(P)\ge \frac14$ for some $u\in V(T)$, or 
there exist  unrelated sets $A, B\subseteq V(T)$ such that   $\omega(A)\ge \frac14$ and  $\omega(B)\ge \frac14$. 
\end{lem}

In the same paper, the authors believe that Lemma~\ref{L}  holds for $\frac13$. This problem has a Ramsey Theory flavor.
In this note, we give an affirmative answer to their question and point out  that  the bound $\frac13$ is tight.

\begin{thm}\label{T1}

Let $T$ be a tree  rooted at $r$ and let $\omega$ be a nonnegative weight function defined on  $V(T)$ with $\omega(T)=1$.  Then there is an $(r, u)$-path $P$ with $\omega(P)\ge \frac13$ for some $u\in V(T)$, or 
there exist  unrelated sets $A, B\subseteq V(T)$ such that   $\omega(A)\ge \frac13$ and  $\omega(B)\ge \frac13$. The bound $\frac13$ is tight.

\end{thm}

\begin{figure}[hbt]
\begin{center}
\includegraphics[scale=1.3]{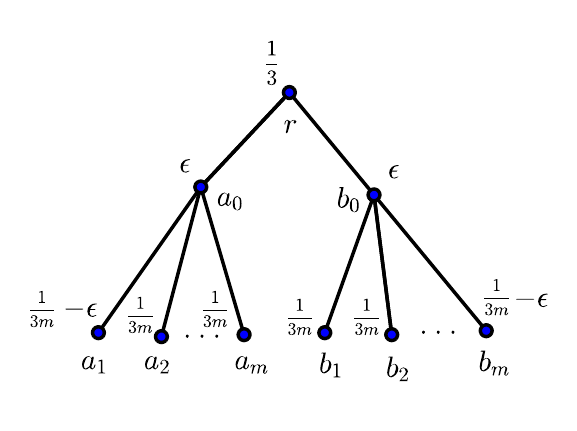}
\caption{Rooted tree $T$.}
\label{T}
\end{center}
\end{figure}

 To see why the bound $\frac13$ is best possible. Let $m\ge2$ be an integer and $\epsilon\ge0$ be a small number with $\epsilon\le\frac1{3m}$.  Let $T$ be  the weighted tree rooted at $r$ as depicted  in Figure~\ref{T}.   Note that  $\omega$ is a positive weight function on $V(T)$  when $\frac1{3m}> \epsilon>0$.  Any path from the root $r$ in $T$ has weight  between $\frac13$ and $\frac13+\frac1{3m}+\epsilon$;  and $T$ has one unique pair of  unrelated sets $A=\{a_0, a_1, a_2, \dots, a_m\}$ and $B=\{b_0, b_1, b_2, \dots, b_m\}$ with $\omega(A)=\omega(B)=\frac13$. The bound $\frac13$ is tight when $m$ is large. \bigskip

\section{Proof of Theorem~\ref{T1}}

 Suppose $T$ has no path from the root $r$ with weight at least $1/3$.  Then $T$ is not a path. Let  $N_G(r)=\{v_1, v_2,\dots, v_s\}$ and $T_1,T_2,\ldots, T_s$ be connected components of $T-r$, where $\omega(T_1)\le \omega(T_2)\le \dots\le \omega(T_s)$.  We call each  $T_i$  a subtree of $T$ rooted at $v_i$ for $1\le i\le s$. 
    And   $T_1, \dots, T_s$  are also called {\it branches}  of $T$ at $r$.    
We next  construct two unrelated sets $A$ and $B$ with desired weights according to the following algorithm:\bigskip

\begin{algorithm}[H]
\SetAlgoLined
\KwData{Vertex weighted tree $T$ with root $r$}
\KwResult{Unrelated sets $A$ and $B$ with desired weights}
Start at the root $r$ with $A=B=\emptyset$  and set $C=\{T_1,T_2,\ldots ,T_s\}$\;
\While{$C\neq \emptyset$}{
\For{$i=1$ to $s-1$}{
Remove $T_i$ from $C$. Add the vertices of $T_i$ to $A$ if $\omega(A)\leq\omega(B)$, and to $B$ otherwise\;
}{
If  $\emptyset\ne\displaystyle\bigcup_{i=1}^{s-1} V(T_i)\subseteq A$ (resp. $\emptyset\ne\displaystyle\bigcup_{i=1}^{s-1} V(T_i)\subseteq B$), color the root $r$ RED (resp. BLUE), otherwise color  the root $r$ GREEN \;
Set $r$ to be the root of $T_s$ and $C$ be the set of connected components of $T_s\backslash r$ with weights sorted in the nondecreasing order\;
}
}
Call the last root $r^*$.  If $\omega(A)\leq \omega(B)$,  add $r^*$ to $A$ and color $r^*$ RED, otherwise add $r^*$ to $B$ and color $r^*$ BLUE. Let $y=r^*$, $x$ be the parent of $y$ and $c$ be the color of $y$\;
\While {$x$ is colored GREEN or c}{
re-color $x$ by the color c if $x$ is colored GREEN and $d_T(x)=2$  \;
Set $y$ to be $x$,  and $x$ be the parent of  $y$ \;
}
Let $u=r^*$\;
\While{$u$ is adjacent to a  vertex $v\not\in A\cup B$  with the same color as $r^*$}{
Add $v$ to $A$  if both $u$ and $v$ are colored RED, and add $v$ to $B$ if both $u$ and $v$ are colored BLUE\;
Set $u$ to  be $v$\;
}

\caption{Building Sets $A$ and $B$}
\end{algorithm} 
\vskip 1cm

 It can be easily checked that  $A$ and $B$ constructed by the above algorithm are unrelated.  Since $T$ is not a path, both  $A$ and $B$ are nonempty.  Let $u$ be the   vertex  in the last step of the algorithm that is added to $A$ or $B$.  According to the algorithm, $u$ is colored RED or BLUE. Let $M$ be the set of all colored vertices of $T$. Then the subgraph $T[M]$ of $T$ induced by $M$ is the unique $(r, r^*)$-path, say $P$, where $r^*$ is the last root as given in the algorithm.  By the algorithm, $T-A\cup B$ is the unique $(r,u^*)$-path, say $P^*$, of $T$, where $u^*$ is the parent of $u$ in $T$.  Clearly, $P^*$ is a subpath of $P$. Let $N=V(P)-V(P^*)$. Then $r^*\in N$ and the vertices of $N$ are all colored by the same color of the root $r^*$. One can see  that if  $u$ is colored RED, then $u\in N\subseteq A$ and the last set of vertices added to $B$ are all uncolored. Similarly, if $u$ is colored BLUE, then $u\in N\subseteq B$ and the last set of vertices added to $A$ are all uncolored.      Since $\omega(P^*)+\omega(A)+\omega(B)=1$ and  $\omega(P^*)<\frac13$, we have\medskip

\noindent (1) \,\, $\omega(A)+\omega(B)>\frac23$.\\

\noindent We next show that $\min\{\omega(A), \omega(B)\}\ge \frac13$.\medskip

Suppose that $\omega(A) \leq \omega(B)$. By (1), $\omega(B)\ge \frac13$.  Assume  $u\in B$.  Then $u$ is colored BLUE and so $r^*$ is also colored BLUE. Thus $N\subseteq B$. Since $r^*$ is added to $B$, we have   $\omega(A)\ge \omega(B-N)$. On the other hand,  $\omega(A)+\omega(B-N)=1-\omega(P)>\frac23$. Thus  $\omega(A)>\frac13$, as desired. So we may assume $u\in A$. Then $u$ is colored RED and so $N\subseteq A$.  Let $D$ be the set of vertices that were last added to $B$. Then $D$ contains only uncolored vertices of $T$. Thus $D=V(Y)$, where  $Y$ is a branch of some  subtree $T^*$ of $T$.  Since $D$ contains only uncolored vertices, by the algorithm, $T^*$ has a branch $X$ with $\omega(Y)\le \omega(X)$ and $X\cap B=\emptyset$.  Let $X^*$ be the set of all vertices that are added to $A$ after the vertices in $D$ were added to $B$. By the algorithm,  $X\subseteq X^*$,  and so $\omega(X^*)\ge\omega(X)\ge \omega(Y)$.  Let $\tilde A=A -X^*$ and $\tilde B=B- Y$.  Since $Y$ is added to $B$, we have $\omega(\tilde A)\ge \omega(\tilde B)$.  Note that $\omega(A) = \omega(\tilde A) + \omega(X^*)$ and $\omega(B) = \omega(\tilde B) + \omega(Y)$. Thus $\omega(B)\le  \omega(\tilde A)+ \omega(X^*)= \omega(A)$.  By (1),   $\omega(A) \ge\frac13$. Hence $\frac13\le \omega(A)\le \omega(B)$, as desired. \\

By a similar argument as above, one can show that  $\min\{ \omega(A), \omega(B)\}\ge \frac 13$ for the case when $\omega(B) \leq\omega(A)$.  This completes the proof of Theorem~\ref{T1}. \hfill\vrule height3pt width6pt depth2pt

\vskip 1cm

\noindent {\bf{Acknowledgement}}: The authors would like to thank the referees for helpful comments.

\end{document}